\documentclass[11pt]{article}
\usepackage[english]{babel}
\usepackage{lineno} 
\usepackage{graphicx,multicol}
\usepackage{epic,eepic,epsfig}
\usepackage{caption}

\usepackage{amsmath,amsfonts,amsthm,amssymb}
\usepackage{pstricks,pstricks-add,pst-math,pst-xkey}

\usepackage{algorithm, algorithmic}

\usepackage{graphicx,pstricks,pst-node,amssymb}

\usepackage{tikz}

\pgfdeclarelayer{edgelayer}
\pgfdeclarelayer{nodelayer}
\pgfsetlayers{edgelayer,nodelayer,main}

\tikzstyle{none}=[inner sep=0pt]
\definecolor{hexcolor0xf81e1c}{rgb}{0.973,0.118,0.110}
\definecolor{hexcolor0x3c00ff}{rgb}{0.235,0.000,1.000}

\tikzstyle{whitevertex}=[circle,fill=white,draw=black, scale = 0.5]
\tikzstyle{redvertex}=[circle,fill=hexcolor0xf81e1c,draw=black, scale = 0.5]
\tikzstyle{bluevertex}=[circle,fill=hexcolor0x3c00ff,draw=black, scale = 0.5]
\tikzstyle{greenvertex}=[circle,fill=green,draw=black, scale=0.5]
\tikzstyle{purplevertex}=[circle,fill=magenta,draw=black, scale=0.5]
\tikzstyle{grayvertex}=[circle,fill=white,draw=gray, scale=0.5]
\tikzstyle{blackvertex}=[circle,fill=black,draw=black, scale=0.5]

\tikzstyle{textbox}=[rectangle,fill=none,draw=none]
\tikzstyle{box}=[rectangle,fill=none,draw=black]

\tikzstyle{arc}=[black, ->]
\tikzstyle{grayarc}=[gray, ->]
\tikzstyle{bluearc}=[blue, ->]
\tikzstyle{grayedge}=[draw=gray]
\tikzstyle{blueedge}=[draw=blue]
\tikzstyle{rededge}=[draw=red]
\tikzstyle{edge}=[draw=black]

\tikzstyle{vertex}=[circle, ,fill=white,draw=black, scale=0.5]

\tikzstyle{10circle}=[circle, scale=10.0,draw=black]
\tikzstyle{10oval}=[ellipse, scale=10.0,draw=black]

\usepackage{amsmath,amsfonts,amsthm,amssymb}
\usepackage{pstricks,pstricks-add,pst-math,pst-xkey}

\usepackage{algorithm, algorithmic}

\usepackage{graphicx,pstricks,pst-node,amssymb}

\usepackage{algorithm, algorithmic}

\begin{document}

\newtheorem{tm}{\hspace{5mm}Theorem}
\newtheorem{prp}[tm]{\hspace{5mm}Proposition}
\newtheorem{dfn}[tm]{\hspace{5mm}Definition}
\newtheorem{lemma}[tm]{\hspace{5mm}Lemma}
\newtheorem{cor}[tm]{\hspace{5mm}Corollary}
\newtheorem{conj}[tm]{\hspace{5mm}Conjecture}
\newtheorem{prob}[tm]{\hspace{5mm}Problem}
\newtheorem{quest}[tm]{\hspace{5mm}Question}
\newtheorem{alg}[tm]{\hspace{5mm}Algorithm}
\newtheorem{sub}[tm]{\hspace{5mm}Algorithm}
\newcommand{\induce}[2]{\mbox{$ #1 \langle #2 \rangle$}}
\newcommand{\2}{\vspace{2mm}}
\newcommand{\dom}{\mbox{$\rightarrow$}}
\newcommand{\ndom}{\mbox{$\not\rightarrow$}}
\newcommand{\compdom}{\mbox{$\Rightarrow$}}
\newcommand{\cdom}{\compdom}
\newcommand{\sdom}{\mbox{$\Rightarrow$}}
\newcommand{\lsd}{locally semicomplete digraph}
\newcommand{\lt}{local tournament}
\newcommand{\la}{\langle}
\newcommand{\ra}{\rangle}
\newcommand{\pf}{{\bf Proof: }}
\newtheorem{claim}{Claim}
\newcommand{\beq}{\begin{equation}}
\newcommand{\eeq}{\end{equation}}
\newcommand{\<}[1]{\left\langle{#1}\right\rangle}

\newcommand{\Z}{\mathbb{$Z$}}
\newcommand{\Q}{\mathbb{$Q$}}
\newcommand{\R}{\mathbb{$R$}}


\title{Gallai-like characterization of strong cocomparability graphs}

\author{Jing Huang\thanks{Department of Mathematics and Statistics,
      University of Victoria, Victoria, B.C., Canada V8W 2Y2; huangj@uvic.ca}}
\date{}

\maketitle

\begin{abstract}
Strong cocomparability graphs are the reflexive graphs whose adjacency matrix can 
be rearranged by a simultaneous row and column permutation to avoid the submatrix 
with rows $01, 10$. Strong cocomparability graphs form a subclass of cocomparability
graphs (i.e., the complements of comparability graphs) and can be recognized in 
polynomial time. 
In his seminal paper, Gallai characterized cocomparability graphs in terms of 
a forbidden structure called asteroids. Gallai proved that cocomparability 
graphs are precisely those reflexive graphs which do not contain asteroids.  

In this paper, we give a characterization of strong cocomparability graphs which
is analogous to Gallai's characterization for cocomparability graphs. We prove that
strong cocomparability graphs are precisely those reflexive graphs which do not 
contain weak edge-asteroids (a weaker version of asteroids). Our characterization 
also leads to a polynomial time recognition algorithm for strong cocomparability 
graphs.  
\end{abstract}

\section{Introduction}

Comparability graphs are a popular and much studied class of graphs 
\cite{gallai,ghouri,gh,golumbic1977,golumbic,gru,hh,mohring}. 
They are the graphs that represent the comparability relation of partial orders 
\cite{trotter} and have been used for the study of optimized compound samples for 
structure-activity correlations of chemical compounds \cite{bkm}. 

Specifically, a graph $G$ is a {\em comparability graph} if it has a transitive 
orientation, that is, the edges of $G$ can be oriented in such a way that for any 
three vertices $x, y, z$ in the resulting oriented graph, if $xy, yz$ are arcs then
$xz$ is also an arc. The complements of comparability graphs are called 
{\em cocomparability graphs}, cf. \cite{golumbic}. 

Cocomparability graphs have a characteristic ordering property, that is, a graph is
a cocomparability graph if and only if it has a vertex ordering $\prec$ such that 
for any three vertices $x \prec y \prec z$, if $xz$ is an edge then at least one of
$xy, yz$ is an edge of the graph. Such a vertex ordering is called a
{\em cocomparability ordering}.

Typically, cocomparability graphs are considered to be {\em reflexive} (i.e., every
vertex is adjacent to itself). Thus the adjacency matrix of a cocomparability graph
has 1's on the main diagonal.
Cocomparability graphs can also be equivalently defined in term of the existence of 
symmetric orderings of their adjacency matrices which do not contain a single 
matrix in specified positions. 

Let $M$ be a symmetric $0,1$-matrix having 1's on the main diagonal. 
A {\em symmetric ordering} of $M$ is a matrix obtained from $M$ by simultaneously 
permuting the rows and columns. Clearly, a symmetric ordering of a symmetric matrix
is again symmetric. Suppose that $M$ is the adjacency matrix of graph $G$. 
Since $M$ has 1's on the main diagonal, $G$ is reflexive. Any vertex ordering of 
$G$ corresponds to a symmetric ordering of $M$, that is, permuting the rows and 
columns of $M$ according to the vertex ordering of $G$. In particular, 
a cocomparability ordering of $G$ corresponds to a symmetric ordering of $M$ which 
does not contain $\left[\begin{matrix}0&1\\
                1&0
                \end{matrix}\right]$
as a submatrix with either 1 entry sitting on the main diagonal. In fact, 
cocomparability graphs are precisely those reflexive graphs whose adjacency 
matrices have such orderings.     

A reflexive graph is called a {\em strong cocomparability} graph if its adjacency
matrix has a symmetric ordering which does not contain
$\left[\begin{matrix}0&1\\
                1&0
                \end{matrix}\right]$ as
a submatrix \cite{hhl}. It follows from the definition that every strong 
cocomparability graph is a cocomparability graph. Both cocomparability graphs and
strong cocomparability graphs are related to interval graphs. According to
Gilmore and Hoffman \cite{gh}, an interval graph if and only if it is both 
a chordal graph and a cocomparability graph. It is proved in \cite{hhl} that 
a graph is an interval graph if and only if it is both a strongly chordal graph and
a strong cocomparability graph. 

Following \cite{hhl,hhlm}, we call the matrix 
$\left[\begin{matrix}0&1\\
                1&0
                \end{matrix}\right]$
the {\em $Slash$ matrix}. A symmetric ordering of $M$ that does not contain 
the $Slash$ matrix as a submatrix is called a {\em symmetric $Slash$-free} 
ordering of $M$. Thus a reflexive graph is a strong cocomparability graph if and 
only if its adjacency matrix has a symmetric $Slash$-free ordering. The vertex 
ordering of a reflexive graph which corresponds to a symmetric $Slash$-free 
ordering of its adjacency matrix is called a {\em strong cocomparability ordering} 
of the graph.  

There is an elegant characterization of cocomparability graphs given by Gallai
\cite{gallai}.  An {\em asteroid} in a graph is a set of 
vertices $x_0, x_1, \dots, x_{2k}$ such that for each $i = 0, 1, \dots, 2k$, 
there is a walk connecting $x_{i+k}$ and $x_{i+k+1}$ which does not contain
any neighbour of $x_{i}$ (subscripts are modulo $2k+1$). 
Gallai \cite{gallai} proved the following:

\begin{tm} \cite{gallai} \label{gallai} 
A reflexive graph is a cocompoarability graph if and only if it does not contain 
an asteroid. 
\qed
\end{tm}

An asteroid consisting of three vertices is called an {\em asteroidal triple}. 
Lekkerkerker and Boland \cite{lb} proved that a graph is an interval graph if and 
only if it is a chordal graph that does not contain an asteroidal triple. 
Graphs which do not contain asteroidal triples are studied in 
\cite{beisegel,cos97,cos99,cs,ducoffe} 

There is also an edge version of asteroids defined in \cite{fhh}. 
An {\em edge-asteroid} in a graph is a set of edges 
$x_0y_0, x_1y_1, \dots, x_{2k}y_{2k}$ such that for each $i = 0, 1, \dots, 2k$, 
there is a walk that begins with the edge $x_{i+k}y_{i+k}$ and ends with the edge 
$x_{i+k+1}y_{i+k+1}$ such that neither $x_{i}$ nor $y_{i}$ is adjacent to a vertex 
in the walk (subscripts are modulo $2k+1$). 
Bipartite graphs which do not contain edge-asteroids are known as
{\em cocomparability bigraphs} and studied in \cite{hhlm}. 
It is shown in \cite{fhh} that a bipartite graph is an 
interval containment bigraph if and only if it is both a chordal bigraph 
and a cocomparability graph, cf. also \cite{huang}.
Note that an asteroid is a special edge-asteroid where each edge is a loop 
(that is, $x_i = y_i$ for each $i  = 0, 1, \dots, 2k$).

In this paper, we introduce the concept of weak edge-asteroids (see Section~\ref{2}
for the definition). Weak edge-asteroids are a weaker version of edge-asteroids in 
the sense that each edge-asteroid (and hence each asteroid) is a weak edge-asteroid.
We will prove the following:

\begin{tm} \label{gallai-like}  
A reflexive graph is a strong cocomparability graph if and only if it does not
contain a weak edge-asteroid. 
\end{tm}

Strong cocomparability graphs can be recognized in polynomial time \cite{hhl}. 
We will show that Theorem \ref{gallai-like} also leads to a polynomial time 
recognition algorithm for strong cocomparability graphs.

\section{Weak edge-asteroids and strong cocomparability graphs} \label{2}

The concept of weak edge-asteroids, to be defined in this section, stems from 
a ``forcing" relation on the set of all ordered pairs of distinct vertices in 
a graph. 

Let $G$ be a reflexive graph with vertex $V(G)$ and edge set $E(G)$. Note that
$E(G)$ includes all {\em loops} $vv$, $v \in V(G)$. Denote by $Z(G)$
the set of ordered pairs $(u,v)$ of distinct vertices of $G$.
For $(u,v), (u',v') \in Z(G)$, we say that $(u,v)$ {\em forces} $(u',v')$,
denoted by $(u,v) \Lambda (u',v')$, if $u = u'$ and $v = v'$ or $uu', vv' \in E(G)$
and $uv', vu' \notin E(G)$. Clearly, $(u,v) \Lambda (u',v')$ if and only if
$(v,u) \Lambda (v',u')$. 

\begin{prp} \cite{hhl} \label{pc}
Let $G$ be a reflexive graph. Then the following statements hold.
\begin{enumerate}
\item If $abc$ is an induced $P_3$ in $\overline{G}$, then $(a,b) \Lambda (c,b)$.
\item If $abcd$ is an induced $P_4$ in $G$, then
      $(a,c) \Lambda (a,d) \Lambda (b,d) \Lambda (a,d) \Lambda (b,c)$.
\item If $abcd$ is an induced $C_4$ in $G$, then $(a,d) \Lambda(b,c)$.
\qed
\end{enumerate}
\end{prp}

We say that $(u,v)$ {\em implies} $(u',v')$, denoted by $(u,v) \sim (u',v')$, if 
there exist walks $u_1u_2 \dots u_k$ and $v_1v_2 \dots v_k$ in $G$ where
$(u_1,v_1) = (u,v)$ and $(u_k,v_k) = (u',v')$ such that
$(u_i,v_i) \Lambda (u_{i+1},v_{i+1})$ for each $i = 1, 2, \dots, k-1$.
It is easy to verify that $\sim$ is an equivalence relation on $Z(G)$.


An {\em invertible pair} in $G$ is a pair of distinct vertices $u, v$ such that
$(u,v) \sim (v,u)$. 

\begin{tm} \cite{hhl} \label{ip}
A reflexive graph is a strong cocomparability graph if and only if it does not
an invertible pair.
\qed
\end{tm}

Let $G$ be a reflexive graph and $uu', vv'$ be edges of $G$ with 
$\{u,u'\} \cap \{v,v'\} = \emptyset$. Note that the edges of $G$ include all loops. We say that $uu'$ {\em avoids} $vv'$ if one of the following
holds:

\begin{itemize}
\item $u = u'$, $v = v'$, and $uv \notin E(G)$;
\item $u = u'$, $v \neq v'$, $uv \notin E(G)$, and $uv' \notin E(G)$;
\item $u \neq u'$, $v = v'$, $uv \notin E(G)$, and $uv' \notin E(G)$; 
\item $u \neq u'$, $v \neq v'$, and $\{u,u',v,v'\}$ induces
      a $2K_2$, a $P_4$, or a $C_4$ in $G$.
\end{itemize}

It is clear from the definition that if $uu'$ avoids $vv'$ then $vv'$ avoids $uu'$.
Also observe that $uu'$ avoiding $vv'$ is equivalent to the property that 
each of $u$ and $u'$ has a non-neighbour in $\{v,v'\}$ and each of $v$ and $v'$ has 
a non-neighbour in $\{u,u'\}$. This observation will be used in the proof of 
Theorem \ref{avoidance}.  

We say that an edge {\em avoids} a walk (with at least one edge) in $G$ if it avoids
every edge of the walk. 

\begin{lemma} \label{align}
Let $G$ be a reflexive graph. Suppose that an edge $vv'$ avoids a walk 
$u_0u_1 \dots u_t$ where $t \geq 1$ in $G$. Then the following statements hold:
\begin{enumerate}
\item If neither $\{u_0,u_1,v,v'\}$ nor $\{u_{t-1},u_t,v,v'\}$ induces a $C_4$,
      then for any $x \in \{v,v'\}$ and $y \in \{v,v'\}$, 
      $(u_0,x) \sim (u_t,y)$.
\item If $\{u_0,u_1,v,v'\}$ induces a $C_4$ and $\{u_{t-1},u_t,v,v'\}$ does not
      induce a $C_4$, then for any $y \in \{v,v'\}$, $(u_0,x) \sim (u_t,y)$
      where $x$ is the unique vertex in $\{v,v'\}$ with $u_0x \in E(G)$.
\item If $\{u_0,u_1,v,v'\}$ and $\{u_{t-1},u_t,v,v'\}$ each induces a $C_4$, 
      then $(u_0,x) \sim (u_t,y)$ where $x$ is the unique vertex in $\{v,v'\}$ 
      with $u_0x \in E(G)$ and $y$ is the unique vertex in $\{v,v'\}$ with 
      $u_ty \in E(G)$.
\end{enumerate}
\end{lemma}
\pf If $v = v'$ (i.e., $vv'$ is a loop), then neither $\{u_0,u_1,v,v'\}$ nor 
$\{u_{t-1},u_t,v,v'\}$ induces a $C_4$, that is, only Statement~1 applies. Since 
$vv'$ avoids every edge of the walk, $vu_i \notin E(G)$ for each $i$. 
By Proposition \ref{pc}, $(u_0,v) \Lambda (u_1,v) \Lambda \cdots \Lambda (u_t,v)$
so $(u_0,v) \sim (u_t,v)$ and Statement~1 holds. Therefore we may assume that 
$v \neq v'$. We prove that the statements hold by induction on $t$.

Consider first the base case $t = 1$. 
Note that $\{u_0,u_1,v,v'\} = \{u_{t-1},u_t,v,v'\}$. If $u_0 = u_1$ (i.e., $u_0u_1$
is a loop), then $\{u_0, u_1, v, v'\}$ does not induce a $C_4$ (and hence only 
Statement~1 applies).  
Since $vv'$ avoids $u_0u_1$, $u_0v, u_0v' \notin E(G)$. By Proposition \ref{pc}, 
$(u_0,v) \sim (u_0,v')$ so Statement~1 holds. If $u_0 \neq u_1$, then 
$\{u_0,u_1,v,v'\}$ induces a $2K_2$, a $P_4$, or a $C_4$. 
In the case when $\{u_0,u_1,v,v'\}$ induces a $2K_2$ or a $P_4$, only Statement~1
applies. It follows from Proposition \ref{pc} that for any $x \in \{v,v'\}$ 
and $y \in \{v,v'\}$, $(u_0,x) \sim (u_1,y)$. Thus Statement~1 holds. On the other 
hand, when $\{u_0,u_1,v,v'\}$ induces a $C_4$, only Statement~3 applies. Either  
$u_0u_1v'v$ or $u_0u_1vv'$ is an induced $C_4$. By Proposition \ref{pc}, 
$(u_0,v) \sim (u_1,v')$ if $u_0u_1v'v$ is an induced $C_4$, and 
$(u_0,v') \sim (u_1,v)$ if $u_0u_1vv'$ is an induced $C_4$. Hence Statement~3
holds. Therefore the statements hold for the base case $t=1$.
Assume now that $t \geq 2$ and the statements hold for any walk of length 
less than $t$. 

Suppose first that neither $\{u_0,u_1,v,v'\}$ nor $\{u_{t-1},u_t,v,v'\}$ induces 
a $C_4$. Let $z$ be the unique vertex in $\{v,v'\}$ with $u_{t-1}z \in E(G)$ if 
$\{u_{t-2},u_{t-1},v,v'\}$ induces a $C_4$; otherwise let $z$ be any vertex in 
$\{v,v'\}$. By the inductive hypothesis, for any $x \in \{v,v'\}$, 
$(u_0,x) \sim (u_{t-1},z)$. Since $\{u_{t-1},u_t,v,v'\}$ does not induce a $C_4$,
by considering the walk $u_{t-1}u_t$ and the inductive hypothesis 
$(u_{t-1},z) \sim (u_t,y)$ for any $y \in \{v,v'\}$. Hence for any 
$x \in \{v,v'\}$ and $y \in \{v,v'\}$, $(u_0,x) \sim (u_t,y)$. 

Suppose next that $\{u_0,u_1,v,v'\}$ induces a $C_4$ and $\{u_{t-1},u_t,v,v'\}$ 
does not induce a $C_4$. Let $x$ be the unique vertex in $\{v,v'\}$ with
$u_0x \in E(G)$. Let $z$ be the unique vertex in $\{v,v'\}$ with 
$u_{t-1}z \in E(G)$ if $\{u_{t-2},u_{t-1},v,v'\}$ induces a $C_4$; otherwise
let $z$ be any vertex in $\{v,v'\}$. By the inductive hypothesis,
$(u_0,x) \sim (u_{t-1},z)$. Since $\{u_{t-1},u_t,v,v'\}$ does not induce a $C_4$, 
$(u_{t-1},z) \sim (u_t,y)$ for any $y \in \{v,v'\}$. Hence for any 
$y \in \{v,v'\}$, $(u_0,x) \sim (u_t,y)$.

Finally, suppose that $\{u_0,u_1,v,v'\}$ and $\{u_{t-1},u_t,v,v'\}$ each induces 
a $C_4$. Let $x$ be the unique vertex in $\{v,v'\}$ with $u_0x \in E(G)$ and 
$y$ be the unique vertex in $\{v,v'\}$ with $u_ty \in E(G)$. Let $z$ be the 
the vertex in $\{v,v'\}$ distinct from $y$. Note that $u_{t-1}z \in E(G)$; 
in particular if $\{u_{t-2},u_{t-1},v,v'\}$ induces a $C_4$ then $z$ is
the unique vertex in $\{v,v'\}$ with $u_{t-1}z \in E(G)$. By inductive hypothesis,
$(u_0,x) \sim (u_{t-1},z)$ and $(u_{t-1},z) \sim (u_t,y)$. Therefore
$(u_0,x) \sim (u_t,y)$. 
\qed

A {\em weak edge-asteroid} in a graph $G$ is a set of edges 
$x_0y_0, x_1y_1, \dots, x_{2k}y_{2k}$ such that for each $i = 0, 1, \dots, 2k$, 
$x_iy_i$ avoids a walk that begins with $x_{i+k}y_{i+k}$ and ends with 
$x_{i+k+1}y_{i+k+1}$ (subscripts are modulo $2k+1$).

\begin{tm} \label{ast}
A reflexive graph $G$ contains an invertible pair if and only if it contains 
a weak edge-asteroid.
\end{tm}
\pf Suppose that $x_0y_0, x_1y_1, \dots, x_{2k}y_{2k}$ form a weak edge-asteroid 
in $G$ (as defined above). 
For each $i = 0, 1, \dots, 2k$, let $u_i \in \{x_i, y_i\}$ and 
$v_{i+k} \in \{x_{i+k},y_{i+k}\}$ be arbitrarily chosen, except when 
$\{x_i, y_i, x_{i+k}, y_{i+k}\}$ induces a $C_4$, $u_iv_{i+k}$ is an edge of 
the induced $C_4$. We may assume without loss of generality that the walk that 
begins with $x_{i+k}y_{i+k}$ and ends with $x_{i+k+1}y_{i+k+1}$ (in the definition
of a weak edge-asteroid) has the first vertex $v_{i+k}$ and the last vertex 
$u_{i+k+1}$. This can be realized by adding $v_{i+k}$ to the begining of the walk
and $u_{i+k+1}$ to the end of the walk if necessary. Thus, by Lemma \ref{align}, 
$(u_i,v_{i+k}) \sim (v_i,u_{i+k+1}) \sim (u_{i+1},v_{i+k+1})$ for each
$i = 0, 1, \dots, 2k$. Hence 
\[(u_0,v_k) \sim (v_0,u_{k+1}) \sim (u_1,v_{k+1}) \sim (v_1,u_{k+2}) \sim
  (u_2,v_{k+2}) \sim \cdots \sim (u_k,v_{2k}) \sim (v_k,u_0)\]
which means that $u_0, v_k$ are an invertible pair in $G$.

Conversely, suppose that $u,v$ are an invertible pair in $G$ certified by
the sequence 
\[(u,v) = (u_0,u_t) \Lambda (u_1,u_{t+1}) \Lambda \cdots \Lambda (u_{t-1},u_{2t-1})
\Lambda (v,u).\]
When $t$ is odd, $u_{2i}u_{2i+1}$ avoids every edge in the walk 
$u_{2i+t-1}u_{2i+t}u_{2i+t+1}u_{2i+t+2}$ for each $0 \leq i \leq t-1$
(subscripts are modulo $2t$).
Hence the edges $u_{2i}u_{2i+1}$, $0 \leq i \leq t-1$, form a weak edge-asteroid
in $G$.
When $t$ is even, $u_iu_{i+1}$ avoids every edge in the walk 
$u_{i+t-1}u_{i+t}u_{i+t+1}$ for each $0 \leq i \leq t-2$,
$u_{t-1}u_t$ avoids every edge in the walk $u_{2t-2}u_{2t-1}u_0u_1$, 
and $u_iu_{i+1}$ avoids every edge in the walk $u_{i-t}u_{i-t+1}u_{i-t+2}$
for each $t \leq i \leq 2t-2$. Hence the edges $u_iu_{i+1}$, $0 \leq i \leq 2t-2$ 
form a weak edge-asteroid in $G$.
\qed

Theorem \ref{gallai-like} now follows immediately from Theorems \ref{ip} and 
\ref{ast}.

As a byproduct, we show that Theorem \ref{gallai-like} leads to a polynomial time
recognition algorithm for strong cocomparability graphs.

Let $G$ be a reflexive graph. Define the {\em avoidance graph} $G^*$ of $G$ as
follows: The vertex set of $G^*$ is $E(G)$ and two vertices $e, f$ of $G^*$ are 
adjacent if $e, f$ avoids each other in $G$. 

\begin{tm} \label{avoidance}
Let $G$ be a reflexive graph. Then $G$ is a strong cocomparability graph if and
only if $G^*$ is a comparability graph.
\end{tm}
\pf In view of Theorems \ref{gallai} and \ref{gallai-like}, it suffices to show
that $G$ has a weak edge-asteroid if and only if the complement $\overline{G^*}$ of 
$G^*$ has an asteroid. 

Suppose first that $G$ has a weak edge-asteroid consisting of edges 
$e_0, e_1, \dots, e_{2k}$. We claim that $e_0, e_1, \dots, e_{2k}$ form an asteroid
in $\overline{G^*}$. Since $e_0, e_1, \dots, e_{2k}$ form a weak edge-asteroid in 
$G$, $e_i$ avoids a walk $W_i$ begins with $e_{i+k}$ and ends with $e_{i+k+1}$ for 
each $i = 0, 1, \dots, 2k$. For any edge $f$ in $W_i$, since $e_i$ avoids $f$ in 
$G$, $e_i$ is not adjacent to $f$ in $\overline{G^*}$. For any two pair consecutive 
edges $f, f'$ in $W_i$, since $f, f'$ share a vertex in $G$, by definition they do 
not avoid and hence they are adjacent in $\overline{G^*}$. Hence the the edges of 
$W_i$ form a walk in $\overline{G^*}$ connecting $e_{i+k}$ and $e_{i+k+1}$ that 
contains no neighbour of $e_i$. Therefore $e_0, e_1, \dots, e_{2k}$ form 
an asteroid in $\overline{G^*}$.

Conversely, suppose that $\overline{G^*}$ has an asteroid consisting of vertices
$e_0, e_1, \dots, e_{2k}$. Denote $e_i = x_iy_i$ for each $i = 0, 1, \dots, 2k$
where $x_i, y_i$ are vertices in $G$. Since $e_0, e_1, \dots, e_{2k}$ form an
asteroid in $\overline{G^*}$, there is a walk $W'_i$ connecting $e_{i+k}$ and 
$e_{i+k+1}$ in $\overline{G^*}$ that contains no neighbour of $e_i$ for each 
$i = 0, 1, \dots, 2k$. We show that $W'_i$ can be modified to a walk in $G$
that begins with $x_{i+k}y_{i+k}$ and ends with $x_{i+k+1}y_{i+k+1}$, and is 
avoided by $x_iy_i$. Consider a pair of consecutive vertices $e, f$ in $W'_i$. 
Denote $e = uu'$ and $f = vv'$ where $u, u', v, v'$ are vertices of $G$.
Since $e$ and $f$ are adjacent in $\overline{G^*}$, they do not avoid each other 
in $G$. Thus some vertex in $\{u,u'\}$ is adjacent to both $v, v'$ or
some vertex in $\{v,v'\}$ is adjacent to both $u, u'$ in $G$. Without loss of 
generality assume that $v$ is adjacent to both $u, u'$ in $G$. We claim that
at least one of $uv, u'v$ avoids $x_iy_i$. Since $e_i = x_iy_i$ is not 
adjacent to $e = uu'$ or $f = vv'$ in $\overline{G^*}$ (by definition of an
asteroid), $x_iy_i$ avoids both $uu'$ and $vv'$ in $G$. Hence none of $u, u', v$ 
is adjacent to both $x_i, y_i$ in $G$. If $x_iy_i$ avoids neither of $uv, u'v$,
then some vertex $a \in \{x_i,y_i\}$ is adjacent to both $u, v$ and 
some vertex $b \in \{x_i,y_i\}$ is adjacent to both $u', v$ in $G$. When $a = b$, 
$a$ is adjacent to both $u, u'$, which contradicts the fact that $x_iy_i$ avoids
$uu'$; when $a \neq b$, $v$ is adjacent to both $x_i, y_i$, which contradicts
the fact $x_iy_i$ avoids $vv'$. Hence $x_iy_i$ must avoids at least one of 
$uv, u'v$ in $G$. That is, there always exists an edge $g$ with one endvertex in
$\{u,u'\}$ and the other endvertex in $\{v,v'\}$ which avoids $x_iy_i$ in $G$.   
To modify the walk $W'_i$, we add such an edge $g$ between any consecutive pair of 
edges $e, f$ in $W'_i$ and if necessary repeat $e$ or $f$ to make sure that
we obtain a walk in $G$. The modified walk still begins with $x_{i+k}y_{i+k}$ and
ends with $x_{i+k+1}y_{i+k+1}$, in which every edge avoids $x_iy_i$. 
This shows that $x_0y_0, x_1y_1, \dots, x_{2k}y_{2k}$ form a weak edge-asteroid
in $G$. 
\qed

In polynomial time one can construct avoidance graphs and check whether they are
comparability graphs \cite{ms}. Thus Theorem \ref{avoidance} implies a polynomial 
recognition algorithm for strong cocomparability graphs.

\begin{cor} \cite{hhl}
Strong cocomparability graphs can be recognized in polynomial time.
\qed
\end{cor}

\section{Further remarks}

For a graph $G$, let $B(G)$ be the bipartite graph with vertex set
$\{v', v'':\ v \in V(G)\}$ and edge set $\{u'v'',\ u''v':\ uv \in E(G)\}$.
Note that when $G$ is reflexive, $v'v''$ is an edge of $B(G)$ for each vertex
$v$ of $G$. 

A {\em trampoline} is a complete graph on $k$ vertices $u_1, u_2, \dots, u_k$ with
$k \geq 3$, together with an independent set of $k$ vertices $v_1, v_2, \dots, v_k$
such that each $v_i$ is adjacent to $u_i$ and $u_{i+1}$ (and to no other vertices).
A chordal graph is {\em strongly chordal} if it does not contain a trampoline as 
an induced subgraph \cite{far}. 

\begin{tm} \cite{far} \label{farber}
A reflexive graph $G$ is a strongly chordal graph if and only if $B(G)$ is
a chordal bigraph.
\qed
\end{tm}

There is a nice comparason in terms of matrix orderings between strongly chordal 
graphs and strong comparability graphs and similarly, between chordal bigraphs and 
cocomparability graphs (see \cite{hhl} for the details).
Strongly chordal graphs are the reflexive graphs whose adjacency matrices admit 
symmetric $\Gamma$-free orderings \cite{far}. 
Chordal bigraphs are the bigraphs whose biadjacency matrices admit $\Gamma$-free 
orderings \cite{af,hsk,lub}, while cocomparability bigraphs are the bigraphs whose 
biadjacency matrices have $Slash$-free orderings \cite{hhlm}. Nevertheless,
a similar equivalence as in Theorem \ref{farber} does not hold for strong
cocomparability graphs. 

\begin{prp} \label{g2b}
If $G$ is a strong cocomparability graph, then $B(G)$ is a cocomparability bigraph.
\end{prp}
\pf Let $M$ be the adjacency matrix of $G$. Then $M$ is the biadjacency matrix of 
$B(G)$. If $G$ is a strong cocomparability graph, then $M$ has a symmetric 
$Slash$-free ordering. Hence $B(G)$ is a cocomparability bigraph.
\qed

The converse of Proposition \ref{g2b} is however not true. For example, $K_{3,3}$ 
is not a strong cocomparability graph as it contains a weak edge-asteroid.
On the other hand, $B(K_{3,3})$ does not contain an edge-asteroid so it is 
a cocomparability bigraph. 

Farber \cite{far} established yet another relationship between strongly chordal 
graphs and chordal bigraphs. For a bigraph $H$, let $H^+$ be the reflexive graph 
obtained from $H$ by completing one colour class of $H$ to a clique and adding 
a loop at each vertex. 

\begin{tm} \cite{far} \label{farb}
A bigraph $H$ is a chordal bigraph if and only if $H^+$ is a strongly chordal 
graph.
\qed
\end{tm} 

The following theorem assembles a similar relationship between strong 
cocomparability graphs and cocomparability bigraphs. 
For a bigraph $H$, let $H^{++}$ be the reflexive graph obtained from $H$ by 
completing both colour classes of $H$ to cliques and adding a loop at each vertex.

\begin{tm}
A bigraph $H$ is a cocomparability bigraph if and only if $H^{++}$ is
a strong cocomparability graph.
\end{tm}
\pf Suppose that $H$ is a cocomparability bigraph. Then the biadjacency matrix of
$H$ has a $Slash$-free ordering $N$. The matrix
$\left[\begin{matrix}{\bf J} & N\\
                {N}^T & {\bf J}
                \end{matrix}\right]$
(where $\bf J$ is an all-ones matrix) is the adjacency matrix of $H^{++}$ 
and does not contain the $Slash$ matrix as a submatrix. Hence $H^{++}$ is 
a strong cocomparability graph.

Conversely, suppose that $H^{++}$ is a strong cocomparability graph. Then the
adjacency matrix has a symmetric $Slash$-free ordering $M$. Deleting the rows
that correspond to the vertices in one colour class and columns that correspond
to the vertices in the other colour class, we obtain the biadjacency matrix 
of $H$ which does not contain the $Slash$ matrix as a submatrix. Hence $H$ is 
a cocomparability bigraph. 
\qed

\end{document}